\documentclass{article}
\usepackage[utf8]{inputenc} % allow utf-8 input
\usepackage[T1]{fontenc}    % use 8-bit T1 fonts
\usepackage{tikz}
\usepackage[style=numeric, sorting=none]{biblatex}
\usepackage{algorithm}
\usepackage{algpseudocode}
\usepackage{float}
\usepackage{amssymb}
\usepackage{amsmath}
\usepackage{hyperref}
\usepackage{mathtools}
\usepackage{graphicx}
\usepackage{subcaption}
\usepackage{siunitx}
\usepackage{booktabs}
\usepackage{tikz-cd}
\usepackage{stix}
\usepackage{bm}
\usepackage{amsbsy}
\usepackage{amsthm}
\usepackage{arxiv}
\usepackage{xspace}
\usepackage{xurl}
\newcommand{\fenics}{\texttt{FEniCs}\xspace}

\newcommand{\velocity}{\bm{v}}

\newcommand\norm[1]{\left\lVert#1\right\rVert}

\newcommand{\inflowBoundary}{\Gamma_{-}}
\newcommand{\characteristicBoundary}{\Gamma_{0}}
\newcommand{\outflowBoundary}{\Gamma_{+}}

\newcommand{\misfit}{{y}}

\newcommand{\normal}{\bm{n}}

\newcommand{\obsO}{\mathcal{B}} % observation operator 
\newcommand{\pto}{\mathcal{F}} % parameter to observable map
\newcommand{\pts}{\mathcal{K}} % parameter to state operator
\newcommand{\som}{\mathcal{M}^+} 
\newcommand{\mta}{\mathcal{Q}} % misfit to adjoint-field operator
\newcommand{\adjoint}{q} % misfit to adjoint-field operator
\newcommand{\ndof}{{n_{\text{dof}}}}
\newcommand{\mdof}{{m_{\text{dof}}}}

\newcommand{\measure}{\mu}

\newcommand{\parameterInitial}{m}

\newcommand{\measureInitial}{\measure}

\newcommand{\ansatzSources}{\mathcal{S}} 
\newcommand{\ansatzSourcesInitial}{\mathcal{S}}

\newcommand{\tobs}{t^{\text{obs}}}

\newcommand{\x}{x}

\newcommand{\ansatzSpace}{\mathcal{V}_h}
\newcommand{\FEMi}{\mathrm{I}}

\newcommand{\Rnn}{\mathbb{R}_{\geq 0}}

\newcommand{\fracSolidus}[2]{#1/#2\,}

\renewcommand{\vec}{\bm}

\newcommand{\tcr}[1]{\textcolor{black}{#1}}

\title{Sparsity-Driven Source Localization in Tomographic Sensing Applications}

\date{} 					% Or removing it

\author{
  \begin{tabular}{@{}c@{\quad}c@{\quad}c@{}}
    \href{https://orcid.org/0009-0009-6325-3578}{\includegraphics[scale=0.06]{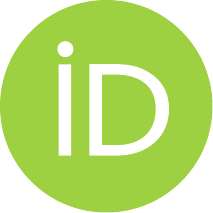}\hspace{1mm}Marco Mattuschka\textsuperscript{$1,$}}\thanks{Corresponding author e-mail: \texttt{marco.mattuschka@dlr.de}} &
    \href{}{\includegraphics[scale=0.06]{orcid.pdf}\hspace{1mm}Noah An der Lan\textsuperscript{$1$}} &
    \href{}{\includegraphics[scale=0.06]{orcid.pdf}\hspace{1mm}Stefanie Schröder\textsuperscript{$2$}} \\
    \noalign{\vskip 3pt}
    \href{https://orcid.org/0000-0002-4819-0048}{\includegraphics[scale=0.06]{orcid.pdf}\hspace{1mm}Arne Ficks\textsuperscript{$2$}} &
    \href{https://orcid.org/0000-0002-2814-0027}{\includegraphics[scale=0.06]{orcid.pdf}\hspace{1mm}Max von Danwitz\textsuperscript{$1$}} &
    \href{https://orcid.org/0000-0002-8820-466X}{\includegraphics[scale=0.06]{orcid.pdf}\hspace{1mm}Alexander Popp\textsuperscript{$1,3$}}
  \end{tabular}
  \\
  \vspace{3pt}
  \textsuperscript{$1$} German Aerospace Center (DLR), Institute for the Protection of Terrestrial Infrastructures,\\
  53757 Sankt Augustin, Germany\\
  \textsuperscript{$2$} Bundeswehr Research Institute for Protective Technologies and CBRN Protection (WIS),\\ 
 29633 Munster, Germany\\
  \textsuperscript{$3$} University of the Bundeswehr Munich, Institute for Mathematics and Computer-Based Simulation (IMCS),\\
  85577 Neubiberg, Germany
}

\renewcommand{\headeright}{}
\renewcommand{\undertitle}{}
\renewcommand{\shorttitle}{Algorithmic Identification of Moving Contaminant Sources}

\begin{document}

\maketitle

\begin{abstract}
\tcr{Hyperspectral standoff detection} systems such as \tcr{Focal Plane Array (FPA)} Fourier Transform Infrared (FTIR) spectrometers provide high spatial resolution in detecting airborne \tcr{chemical} contaminants that are invisible to the human eye but potentially hazardous\cite{10.1117/12.2518834}. When two such systems are operated simultaneously with a suitable opening angle, they enable tomographic reconstruction of contaminant plumes with improved spatial and temporal accuracy\cite{spie:133423fc396d96d6105c143c6a891646b5384334}.
This work presents a mathematical model of these measurement capabilities and an algorithm to identify, localize, and quantify contaminant release sources. The objective is to develop a a tool that reconstructs release locations and predict the future plume evolution from standoff measurement data, thereby supporting early warning and situational awareness in hazardous material release scenarios.
The transport of contaminants is modeled by an advection-diffusion equation, and the corresponding inverse problem for source identification is formulated accordingly. Owing to the severe ill-posedness and underdetermination of the problem, a sparsity-promoting regularization approach is employed together with a high-performance optimization algorithm\cite{MATTUSCHKA2026118854}. To incorporate the tomographic measurement data into the discrete formulation, a level-set description of a threshold concentration 
%immersed boundary method\cite{PESKIN1972252} 
is used, allowing the measurements to be represented independently of the computational mesh and avoiding costly remeshing procedures.
\end{abstract}


\section{Introduction}

\begin{figure}
\centering
\begin{subfigure}{0.32\textwidth}
\centering
\includegraphics[width=0.8\linewidth]{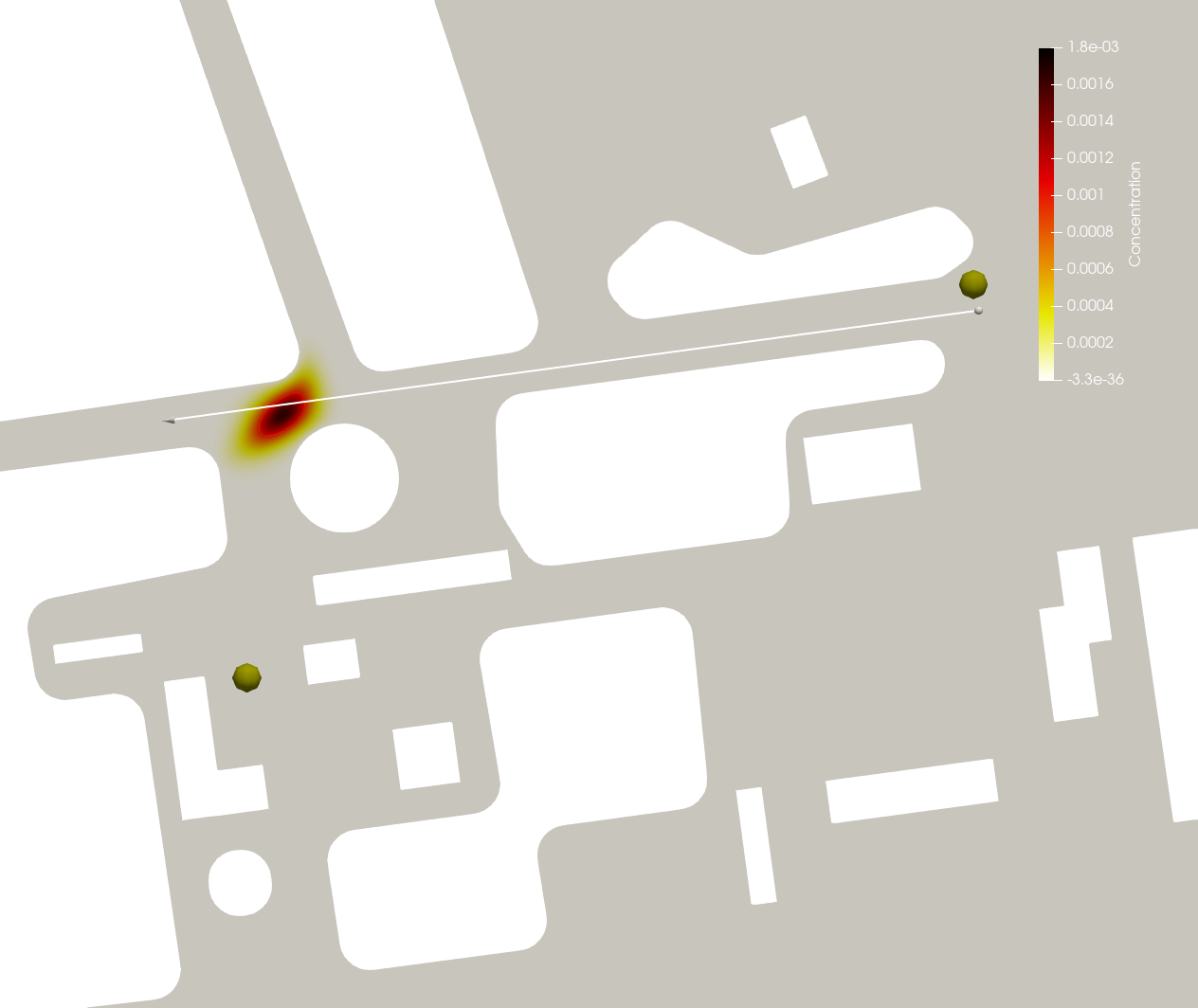}
\end{subfigure}
\begin{subfigure}{0.32\textwidth}
\centering
\includegraphics[width=0.8\linewidth]{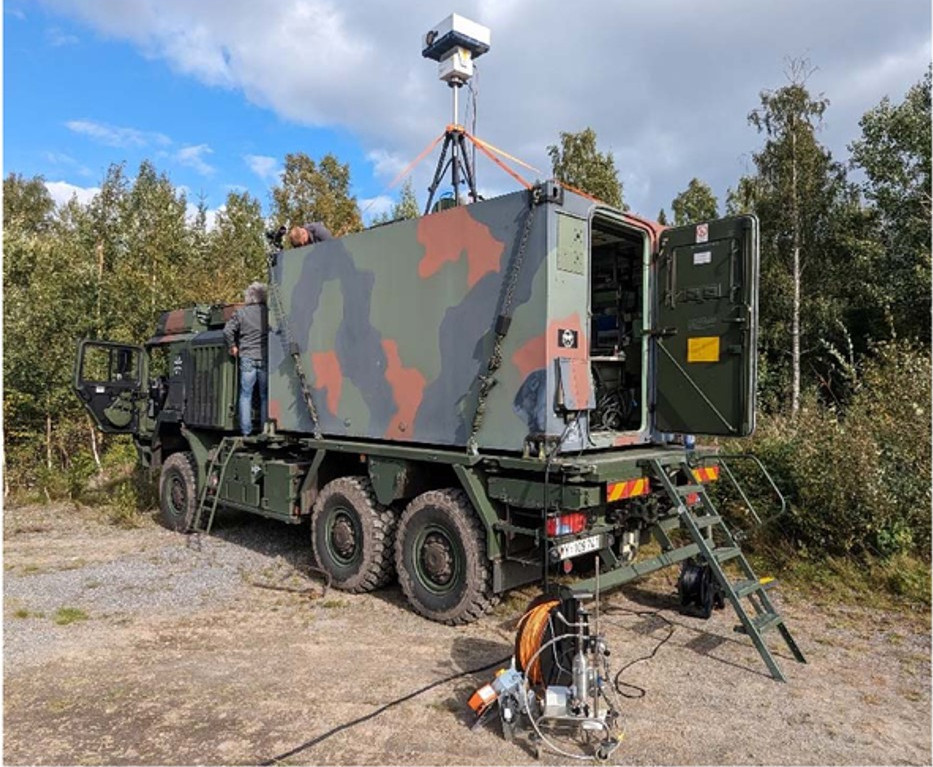}
\end{subfigure}
\begin{subfigure}{0.32\textwidth}
\centering
\includegraphics[width=.8\linewidth]{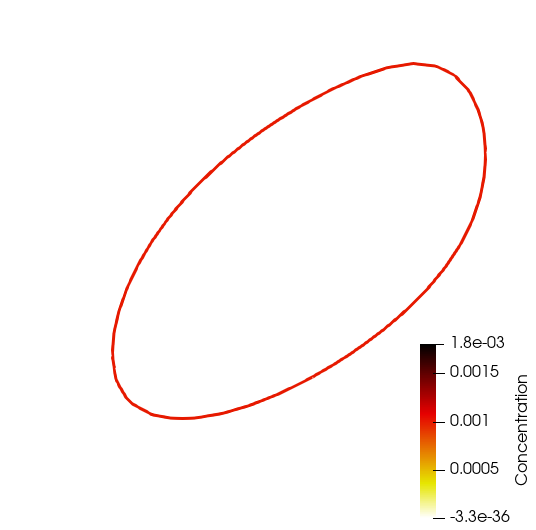}
\end{subfigure}
\caption{Modeled measurement setup with two \tcr{FPA-}FTIR measuring stations (\tcr{green-yellow dots}, left), \tcr{FPA-FTIR hyperspectral imaging measurement device in another} campaign \cite{spie:133423fc396d96d6105c143c6a891646b5384334}(top middle), synthetic signal (extracted contour) representing gas tomography result in the modeled setup (right)}
\label{fig:hi90}
\end{figure}

%\subsection{Motivation}
The threat posed by airborne \tcr{chemical hazards, such as chemical warfare agents (CWAs) or toxic industrial chemicals (TICs)}, is a growing concern in both civilian and military contexts. \tcr{Substances} releases may occur accidentally through industrial leaks or spills, or may be caused intentionally by acts of sabotage or terrorism. The resulting \tcr{highly} toxic gas clouds may be invisible to the human eye, yet can rapidly expand to threaten large areas. Releases of hazardous \tcr{chemical substances} can inflict mass casualties and create widespread disruption in targeted societies. The rapid and accurate identification of the contaminant source and spread is crucial for an effective response and mitigation of the threat. In emergency situations, decision-makers require reliable and timely information about the current state of contamination in order to initiate appropriate countermeasures\cite{ZHANG2025118474}. Situational awareness can be achieved by combining sensor data with physics-based models and suitable inference algorithms. 

%\subsection{Literature Review}

%\subsubsection{Review on sensors}
The \tcr{detection, identification and monitoring} of hazardous airborne \tcr{chemicals} has been an active research area for decades, primarily driven by \tcr{Chemical, Biological, Radiological and Nuclear} (CBRN) defense requirements\cite{Sferopoulos}. Standoff hyperspectral imaging detectors are particularly suited for monitoring a large area within their line of sight\cite{Harig2001}, while remotely controlled point detectors may be used to survey specific areas \tcr{of strategic interest from a safe distance.}\cite{spie:133423fc396d96d6105c143c6a891646b5384334,10168106}. A favorable spatial and temporal detection resolution is achieved \tcr{by using FPA-FTIR-spectrometers}\cite{10.1117/12.692922}. Early work on passive standoff detection technologies established the viability of \tcr{FPA-}FTIR-based sensing for toxic industrial chemicals and chemical warfare agents at operationally relevant distances\cite{10.1117/12.2518834}. Operating two \tcr{FPA-FTIR-based hyperspectral imaging systems} in a tomography setup enables both early warning and a spatial reconstruction of the cloud in the field of view (see \autoref{fig:hi90}(top left) and  \cite{10.1117/12.2663155,spie:133423fc396d96d6105c143c6a891646b5384334} for further details of the measurement system operation).

%\subsubsection{Review on Algorithms}
Physics-based models and suitable inference algorithms have the potential to extend the mapping of a threat zone beyond the field of view of the sensor systems. Various approaches for solving the inverse problem of identifying an initial contaminant release from (partial) measurements of the current state have been proposed in the literature. A Bayesian Markov chain Monte Carlo (MCMC) framework is presented in\cite{Albani.2021}. From a Bayesian perspective, challenges arise in the computation of maximum-a-posteriori estimators, where the specific choice of data misfit and regularization terms is closely related to the assumed noise model and the prior distribution of the unknown source, respectively\cite{Stuart.2010,Nitzler.2022}. In contrast, we formulate the source identification problem as a PDE-constrained optimization problem, where the advection-diffusion equation acts as the governing constraint and the scarcity of measurement data are compensated by suitable regularization techniques\cite{Gorelick.1983,Tsai.2014,Villa.2021}. In this context, Tikhonov regularization with weighted $L^2$-type penalty terms leads to quadratic minimization problems that can be solved efficiently using Newton-type methods\cite{Wiedemann.2024,Petra.2011,Villa.2021,Wu.2023,Alexanderian.2014, Alexanderian.2018, Attia.2018,Wogrin.2023}.

However, in realistic scenarios, contaminant sources are often highly localized, i.e., concentrated around a finite number of release locations, while the available measurements are typically sparse in space, for example in the form of time series recorded by a small number of fixed \tcr{point} sensors. Due to the smoothing properties of the advection-diffusion equation, weighted $L^2$-type regularization is generally not well suited to capture such sparse structures.

To address this issue, we model the unknown contaminant source as a superposition of finitely many atoms of prescribed shape and subsequently relax this representation to parameters described by integrals of a shape function with respect to an unknown positive Radon measure. The resulting identification problem is treated within a variational regularization framework, leading to a convex optimization problem. For its numerical solution, we employ the Primal-Dual Active Point (PDAP) strategy introduced in\cite{Pieper.2021,MATTUSCHKA2026118854}, which efficiently approximates sparse minimizers by alternating between greedy updates of source locations and optimization of the corresponding source intensities.

%\subsection{Scientific Novelty}

%This work adapts the recently presented algorithmic framework for sparsity-driven source localization \cite{MATTUSCHKA2026118854} to combination of sensor data and algorithm. 

The present work addresses a gap at the intersection of two established research lines. While sparse regularization methods for advection-diffusion source identification have been developed and validated on synthetic point sensor configurations, and while standoff \tcr{FPA-FTIR} tomography has been demonstrated as a capable sensing architecture for airborne contaminant detection in field conditions, no prior work has combined these two paradigms into a unified framework for real-time source localization. Specifically, the scientific novelty of this work is twofold. First, we adapt the Primal-Dual-Active-Point strategy of \cite{MATTUSCHKA2026118854} to the tomographic measurement setup of dual FPA-FTIR standoff sensing systems. Therefore, we replace the point sensor observation model with a level-set description of a threshold concentration which corresponds to the information that can be obtained from the measurement system. Second, we demonstrate numerically that the resulting sparse inverse problem remains tractable. Together, these contributions establish sparsity-driven tomographic source localization as a computationally feasible approach.

%\subsection{Paper Organization}
The remainder of this paper is organized as follows. \autoref{sec:math} provides a mathematical framework for the problem, outlining the  equations that we use to model the behavior of airborne hazardous substances. In \autoref{sec:method}, we specialize the Primal-Dual Active Point (PDAP) algorithm for source identification based on level-set information. Numerical results obtained with the proposed method are presented in \autoref{sec:result}, where we demonstrate the feasibility of the approach in a numerical experiment. Finally, in \autoref{sec:outlook}, we summarize the main findings and provide an outlook on further developments and field trials to validate the method.

\section{Mathematical Framework}\label{sec:math}
\begin{figure}
\centering
\begin{subfigure}{0.48\textwidth}
\centering
\includegraphics[width=1.0\linewidth]{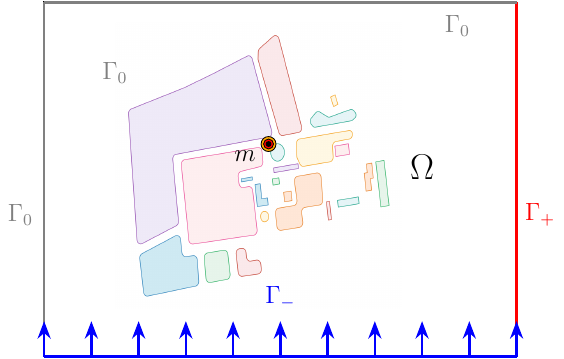}
\end{subfigure}
\begin{subfigure}{0.42\textwidth}
\centering
\includegraphics[width=.9\linewidth]{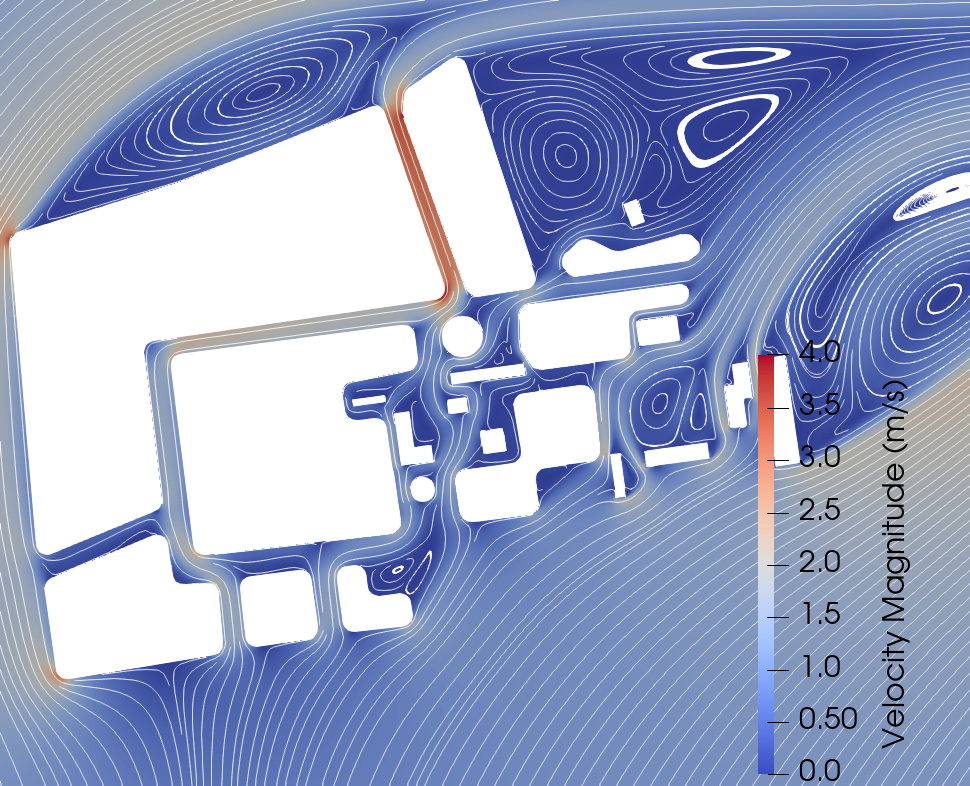}
\end{subfigure}
\caption{Shows the test area $\Omega$ (left) and the estimated wind field (right) }
\label{fig:domain}
\end{figure}
\begin{figure}
\centering
\begin{subfigure}{0.48\textwidth}
\centering
\includegraphics[width=0.80\linewidth]{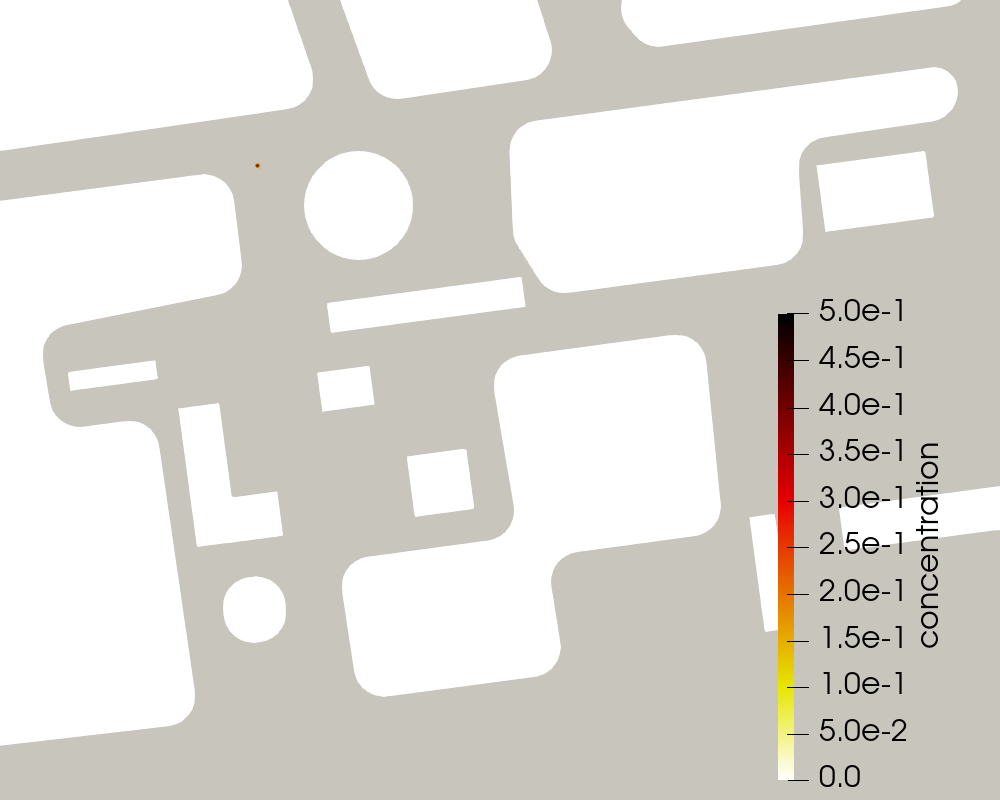}
\end{subfigure}
\begin{subfigure}{0.48\textwidth}
\centering
\includegraphics[width=0.80\linewidth]{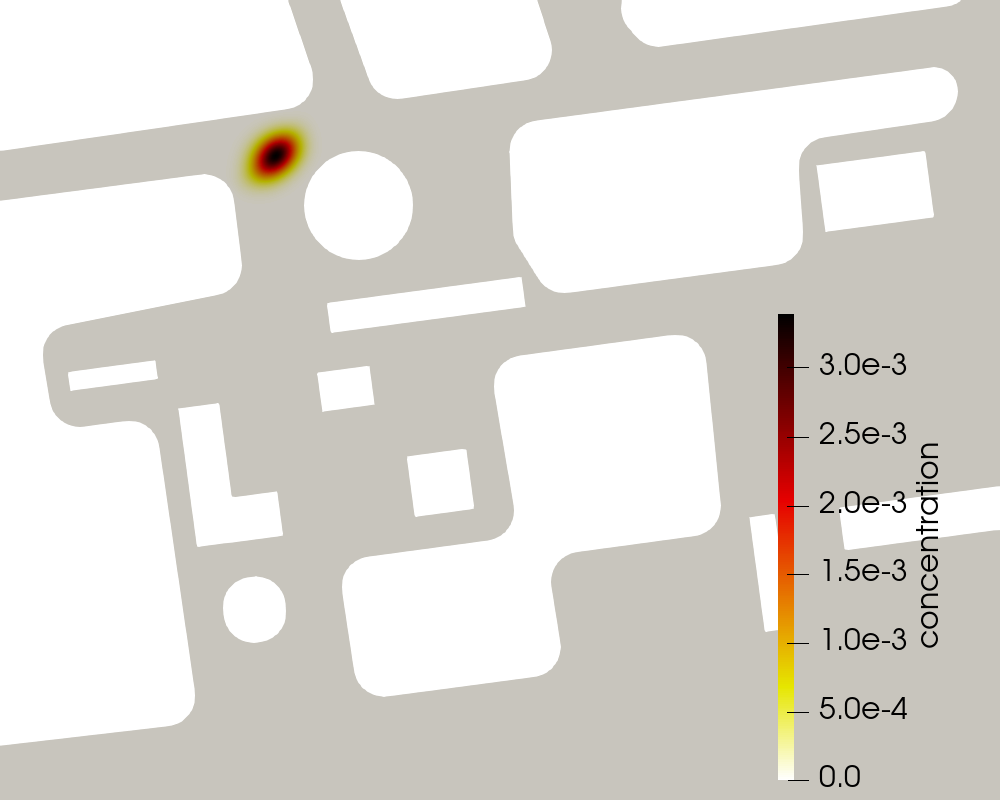}
\end{subfigure}
\begin{subfigure}{0.48\textwidth}
\centering
\includegraphics[width=0.80\linewidth]{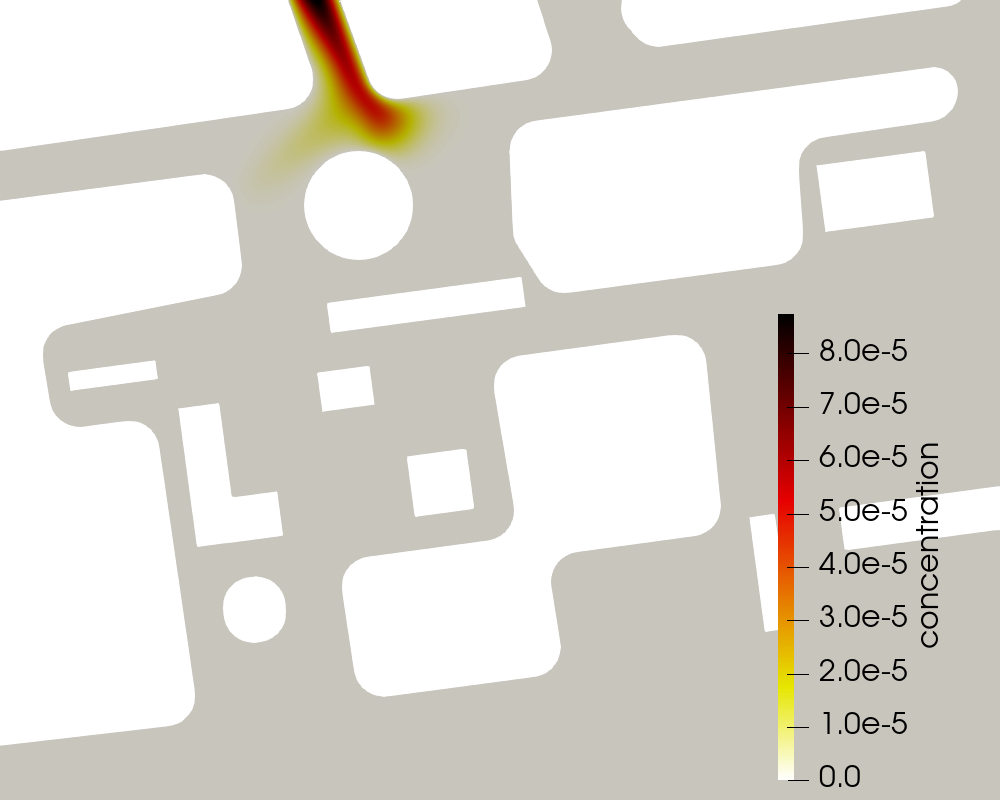}
\end{subfigure}
\begin{subfigure}{0.48\textwidth}
\centering
\includegraphics[width=0.80\linewidth]{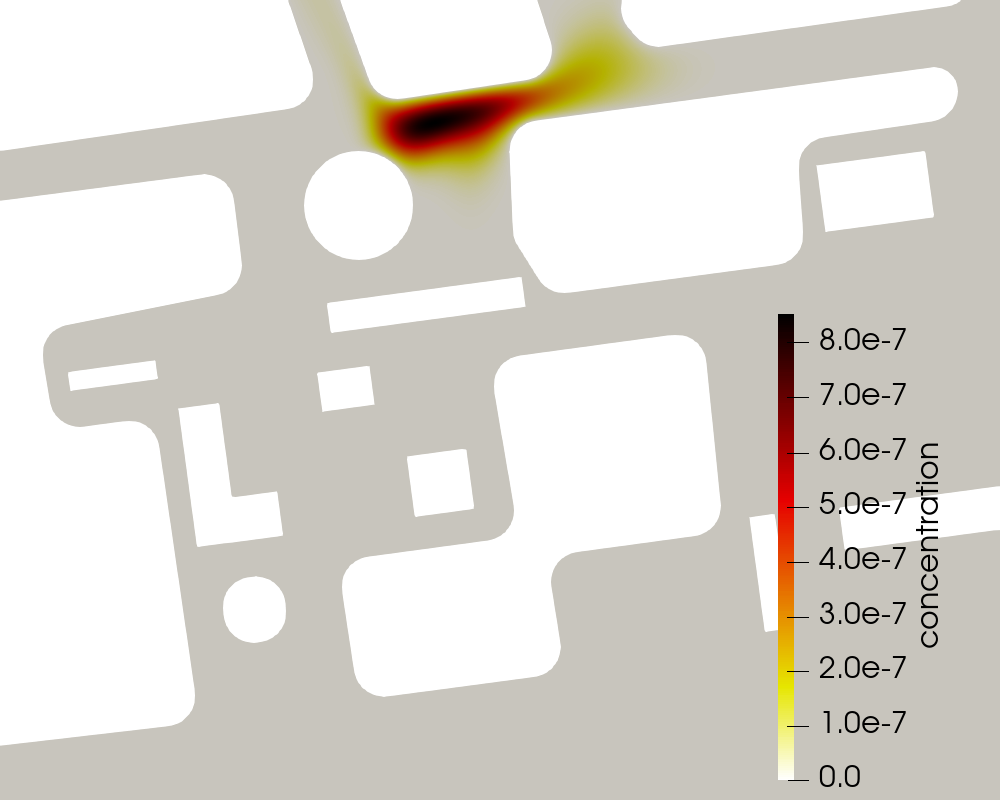}
\end{subfigure}
\caption{Numerical simulation of contamination dispersion at the beginning of the release ($t=\SI{0.0}{\s}$, top left), at the measurement time ($t=\SI{5.0}{\s}$, top right) , after $t=\SI{40.0}{\s}$ (bottom left)  and at the final simulation time ($t=\SI{100.0}{\s}$, bottom right).}
\label{fig:fwd}
\end{figure}
The underlying transport process is described by the advection-diffusion equation with diffusion coefficient $\kappa$ and wind velocity field $\velocity$. The boundary $\partial\Omega$ is decomposed into the inflow boundary $\inflowBoundary$, the free-slip boundary $\characteristicBoundary$, and the outflow boundary $\outflowBoundary$. These boundaries are characterized by the sign of the normal component of the velocity field:
\[
\velocity \cdot \normal < 0 \quad \text{on } \inflowBoundary, \qquad
\velocity \cdot \normal = 0 \quad \text{on } \characteristicBoundary, \qquad
\velocity \cdot \normal > 0 \quad \text{on } \outflowBoundary,
\]
where $\normal$ denotes the outward-pointing unit normal vector on $\partial\Omega$, cf.\cite{Elman.2020}. A schematic illustration of this setting is shown in \autoref{fig:domain}. The resulting advection-diffusion problem with corresponding boundary conditions is given by
\begin{equation}\label{eq:forward_equation}
\begin{aligned}
u_t-\kappa\Delta u + \velocity\cdot\nabla u &= 0 &\qquad&\text{in}\ (0,T)\times\Omega,\\
  \kappa\nabla u \cdot \normal &= 0 &&\text{in}\ (0,T)\times \outflowBoundary ,\\
  u&= 0 &&\text{in}\ (0,T)\times (\inflowBoundary \cup \characteristicBoundary),\\
  u(0,\cdot) &= \parameterInitial &&\text{in}\ \Omega.
\end{aligned}
\tag{$\mathcal{P}_{\pto}$}
\end{equation}
The abstract space of admissible parameters will be denoted by $\mathrm{D}$ and will be specified later. An illustration of the corresponding forward simulation is shown in \autoref{fig:fwd}. 

In the following, we present the extension of the approach introduced in \cite{MATTUSCHKA2026118854}. To this end, for a given real physical contaminant field $u_{\text{true}}$, measured at a height of $\SI{2}{\metre}$ above the ground, \autoref{eq:forward_equation} is assumed to be satisfied, and the \tcr{FPA-FTIR} hyperspectral imaging device, cf. \autoref{fig:hi90}\tcr{(top right)}, is able to determine a level set at concentration value $c$ at time $\tobs \in (0,T]$, i.e.,
\[
\Gamma_{\tobs}[u_{\text{true}}]
    := \{ x \in \Omega \mid u_{\text{true}}({\tobs},x) = c \}.
\]
For simplicity of notation, we write $\Gamma_{\tobs} := \Gamma_{\tobs}[u_{\text{true}}]$, although the set depends on the unknown true solution $u_{\text{true}}$ and assume that $\Gamma_{\tobs}$ is smooth embedded. Neglecting sensor noise, the observation operator at time ${\tobs}$ is defined by $\obsO_{\tobs}(u) := u({\tobs},\cdot)\big|_{\Gamma_{\tobs}}$, which is well-defined due to the regularity of the solution. Our data assimilation step consists of reconstructing the initial condition $\parameterInitial$ such that the simulated concentration matches the measured concentration values along the observed level sets. This leads to the minimization problem
\begin{equation}
\label{eq:objectiv}
\begin{aligned}
    \min_{\parameterInitial \,\text{ admissible}}
    \;
    \int_{\Gamma_{\tobs}}
    \left(
        u(\tobs,z) - c
    \right)^2
    \, \mathrm{d}s(z),
\end{aligned}
\end{equation}
subject to $u$ solving \autoref{eq:forward_equation}. Accordingly, we define the linear and continuous parameter-to-state operator 
\[\pts \colon \mathrm{D} \to C^0([\tobs,T]\times \bar{\Omega}),\qquad \text{defined by }\pts(m)=u,\]
where $u$ solves \autoref{eq:forward_equation}. The corresponding parameter-to-observable operator is given by
\[
\pto \colon \mathrm{D} \to C^0(\Gamma_{\tobs}),
\qquad
\pto(m) = \obsO_{\tobs} \circ \pts(m).
\]
Given a smooth \textit{misfit} $\misfit \in C^0(\Gamma_{\tobs})$, e.g., $\misfit = \obsO_{\tobs}(u) - c$, the associated \textit{misfit-to-adjoint} map is given by $\mta(\misfit) = \adjoint$, where $\adjoint$ is the solution to the final value problem
\begin{equation}\label{eq:adjoint_equation}
\begin{aligned}
  -\adjoint_t - \kappa \Delta \adjoint - \operatorname{div}(\adjoint \velocity) &=  \, \tcr{\delta_{\Gamma_{\tobs}}[\misfit]}
  && \text{in } (0, T) \times \Omega, \\
  (\velocity \adjoint + \kappa \nabla \adjoint) \cdot \normal &= 0 
  && \text{on } (0, T) \times (\inflowBoundary \cup \characteristicBoundary), \\
  \adjoint &= 0 
  && \text{on } (0, T) \times \outflowBoundary, \\
  \adjoint(T, \cdot) &= 0 
  && \text{in } \Omega. 
\end{aligned}
\tag{$\mathcal{P}_{\mta}$}
\end{equation}
The line delta distribution can be expressed informal as
\[
\delta_{\Gamma_{\tobs}}[\misfit](x)
    :=
    \int_{\Gamma_{\tobs}}
    \misfit(z)\,\delta(x-z)\, \mathrm{d}s(z) \qquad \text{for } x \in \Omega,
\]
where $\delta$ denotes the Dirac delta distribution, cf. \cite{PESKIN1972252}. More precisely, $\delta_{\Gamma_{\tobs}}[\misfit]$ is defined by the identity
\[
\int_{\Omega}\delta_{\Gamma_{\tobs}}[\misfit](x)\, f(x)\, \mathrm{d}x=\int_{\Gamma_{\tobs}}\misfit(z) f(z)\, \mathrm{d}s(z) \qquad
\text{for every test function } f \in C_c^\infty(\Omega).
\]

\section{Method: Source identification based on level-set information}\label{sec:method}
\begin{figure}
\centering
\begin{subfigure}{0.48\textwidth}
\includegraphics[width=.9\linewidth]{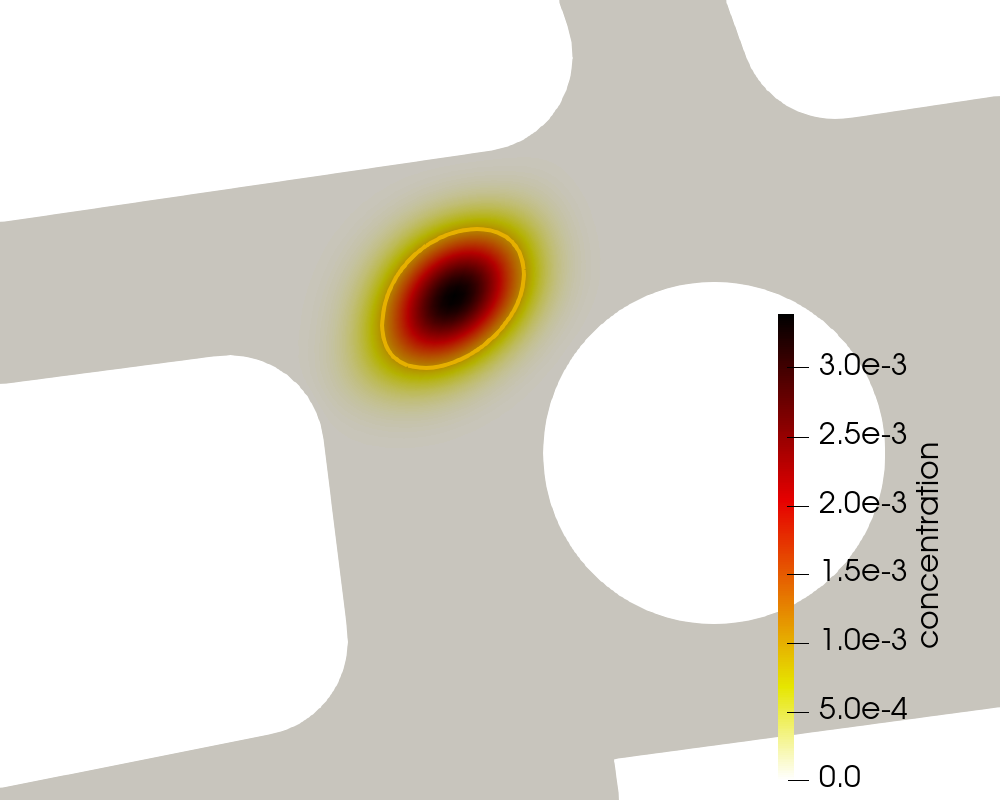}
\end{subfigure}
\begin{subfigure}{0.48\textwidth}
\includegraphics[width=.9\linewidth]{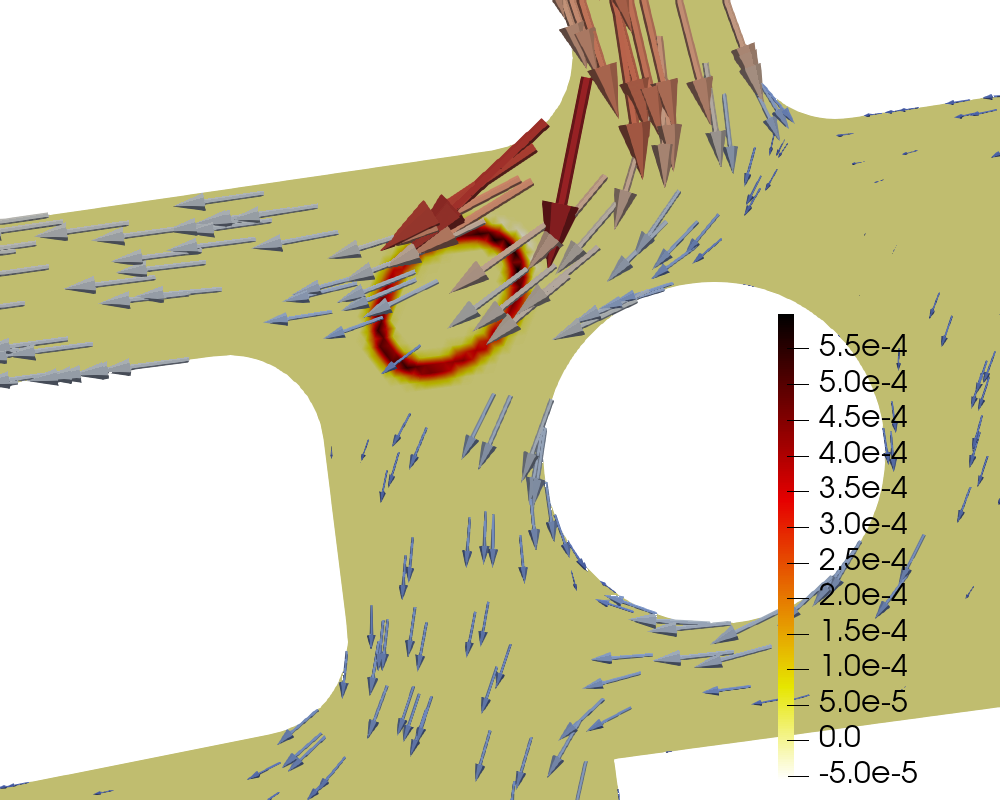}
\end{subfigure}
\caption{Shows the contaminant after $\SI{5.0}{\s}$ seconds and the contour ${\Gamma_{\SI{5.0}{\s}}}=\{x \in \Omega | u_\text{true}(t,x)= 10^{-3}\}$ captured by the measurement system  (yellow ellipse, left) and \tcr{the solution q of the adjoint problem (\ref{eq:adjoint_equation}) for $c=10^{-3}$ at $t=\tobs$, computed using the inverted wind field (right).}}
\label{fig:contur}
\end{figure}
If we consider \autoref{eq:objectiv}, the resulting inverse problem is highly underdetermined, since the initial condition is represented by a sufficiently regular parameter field $\parameterInitial$. To reduce the dimensionality of the problem, we follow the approach proposed in\cite{Pieper.2021,MATTUSCHKA2026118854} and assume that the qualitative structure of the initial condition is known a priori, while only the source locations and amplitudes remain unknown. More precisely, for a prescribed shape function $\ansatzSourcesInitial : \Omega \times \Omega \to \mathbb{R}$, we consider the optimization problem
\begin{equation*}
\begin{aligned}
    \min_{x \in \Omega,\;\lambda \in \Rnn} \frac{1}{2}\,
    \int_{\Gamma_{\tobs}}
    \left(
        \pto(\lambda\ansatzSourcesInitial(x,\cdot))
        (\tobs,z)
        - c
    \right)^2
    \, \mathrm{d}s(z) ,
\end{aligned}
\end{equation*}
which results in a highly nonlinear parameter identification problem. To overcome this difficulty, we assume that the initial condition can be represented as a conic combination of Dirac delta functionals,
\begin{equation}
\label{eq:finiteSum}
\parameterInitial =\int_{\bar{\Omega}} \ansatzSourcesInitial(z,\cdot) \, \mathrm{d}\measureInitial(z), \qquad \text{where } \measureInitial=\sum_{i=1}^{N} \lambda_i\, \delta_{\x_i^{}}.
\end{equation}
We then extend this representation to the space of Radon measures by defining
\begin{equation}
\label{eq:linearparammaps}
\parameterInitial(\measureInitial)= \int_{\bar{\Omega}} \ansatzSourcesInitial(z,\cdot) \, \mathrm{d}\measureInitial(z) \text{ and }\hat{\pto}(\measureInitial):=\pto(\parameterInitial(\measureInitial))  \text{ for a positive Radon measure } \mu .
\end{equation}
Then, the inverse problem can be formulated as the following well-posed optimization problem over the space of positive Radon measures $\som(\bar{\Omega})$:
\begin{equation}
\label{eq:sparseObjective}
\min_{\measureInitial \in \som(\bar{\Omega})} J(\measureInitial):= \frac{1}{2}\,\|\hat{\pto}(\measureInitial)-c\|_{L^2(\Gamma_{\tobs})}^2+\alpha \,\measureInitial(\bar{\Omega}), \tag{$\mathcal{P}_{\mathcal{M}}$}
\end{equation}
where $\alpha > 0$ denotes a regularization parameter.
For the existence and regularity of a minimizer, the SUPG-stabilized finite element discretization,  as well as the numerical determination by the PDAP algorithm, we refer to \cite{MATTUSCHKA2026118854,Danwitz.2023,Brooks.1982} and restrict ourselves here to the discretization of the adjoint problem (\ref{eq:adjoint_equation}) or, more specifically, to the discretization of the distribution $\delta_{\Gamma_{\tobs}}[\misfit]$, which appears as source term on the right-hand side of the equation. We use continuous Lagrange nodal basis functions defined by
$
\ansatzSpace = \operatorname{span}\{\phi_{{1}}, \dots, \phi_{\ndof}\}
$
associated with the nodes $\{p_{{1}}, \dots, p_{\ndof}\}$. 
The mass matrix in $\Omega$ is given by
\[
M^{\Omega}_{ij} := \int_{\Omega} \phi_i(x)\,\phi_j(x)\,\mathrm{d}\Omega(x).
\]
An identification between functions in \(L^2(\Omega)\) and their finite element coefficient vectors is required. For this purpose, let the map $\FEMi: \mathbb{R}^\ndof \rightarrow \ansatzSpace$ be defined by $\FEMi(\vec{u}) = \sum_{i=1}^\ndof u_i \phi_i$. Together with the scalar product $\langle \vec{a}, \vec{b} \rangle_{M^{\Omega}} := \langle M^{\Omega} \vec{a}, \vec{b} \rangle_{\mathbb{R}^\ndof}$, the map $\FEMi$ becomes an isometry, i.e., 
$ \langle \FEMi(\vec{u}), \FEMi(\vec{v}) \rangle_{L^2(\Omega)} = \langle \vec{u}, \vec{v} \rangle_{M^{\Omega}}$. Next, consider a discrete representation of the embedded curve \[
\Gamma_{\tobs}=\{\gamma_1,\dots,\gamma_{\mdof}\},
\]
together with a continuous piecewise linear nodal basis $\{\psi_1,\dots,\psi_{\mdof}\}$ defined on the corresponding one-dimensional mesh of the curve. Let $M^{\Gamma_{\tobs}}\in\mathbb{R}^{\mdof\times\mdof}$ denote the associated mass matrix on \(\Gamma_{\tobs}\). For simplicity and improved numerical robustness, each point \(\gamma_j\) \tcr{is assumed to lie at the barycenter $\gamma_j\approx\frac13(p_{j_0}+p_{j_1}+p_{j_2})$ of a finite element  $K_j=\operatorname{conv}\{p_{j_0},p_{j_1},p_{j_2}\}\subset\Omega$.} Since linear finite elements are employed, the basis functions satisfy
\begin{equation}
C_{ji}
:=\phi_i(\gamma_j)
=
\begin{cases}
\frac13,
& i\in\{j_0,j_1,j_2\},
\\[0.3em]
0,
& \text{otherwise},
\end{cases}
\end{equation}
which defines the sparse trace matrix $C \in \mathbb{R}^{\mdof\times\ndof}$.
The discrete trace operator is therefore given by $\obsO_{\tobs}(\FEMi(\vec u))=\sum_{j=1}^{\mdof}\sum_{i=1}^{\ndof} C_{ji}u_i\psi_j$. Consequently, the discrete representation of the singular source term appearing on the right-hand side of the adjoint problem \eqref{eq:adjoint_equation} is obtained from
\[
\left\langle
\delta_{\Gamma_{\tobs}}[\FEMi(\vec{\misfit})],
\phi_l
\right\rangle
=
\sum_{k=1}^{\mdof}
\int_{\Gamma_{\tobs}}
\psi_k\,
\obsO_{\tobs}(\phi_l)
\,\mathrm ds\,
\misfit_k .
\]
Using the matrix representation of the trace operator, this becomes
\[
\sum_{j=1}^{\mdof}\sum_{k=1}^{\mdof}C_{jl}\,M^{\Gamma_{\tobs}}_{jk}\,\misfit_k=\left(C^\top M^{\Gamma_{\tobs}}\vec{\misfit}\right)_l.
\]
Hence, the discrete contour source term is represented by $C^\top M^{\Gamma_{\tobs}}\vec{\misfit}$, whose spatial distribution is illustrated in \autoref{fig:contur}.

\section{Numerical Results}\label{sec:result}
\begin{figure}
\begin{subfigure}{0.48\textwidth}
\includegraphics[width=1.0\linewidth]{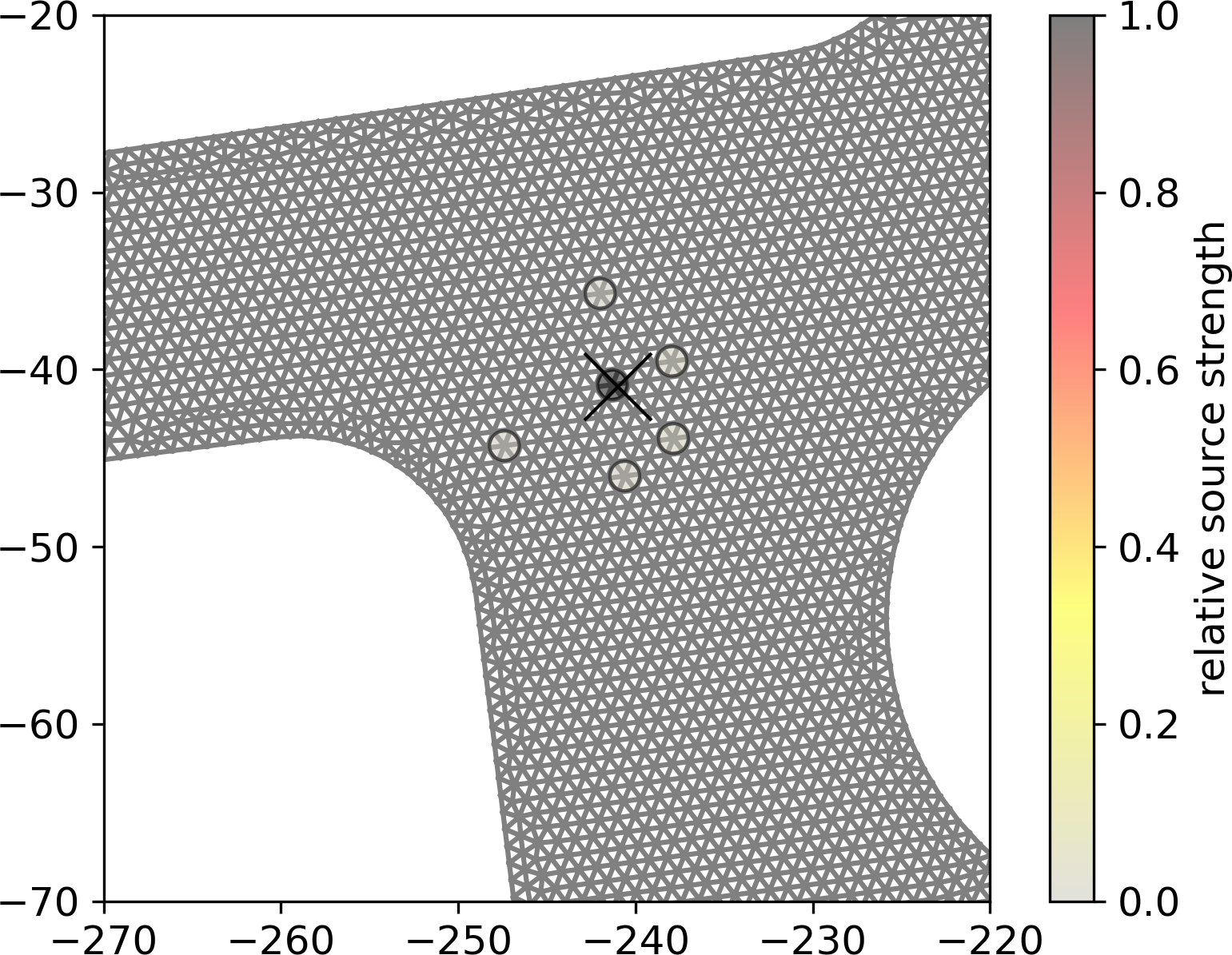}
\end{subfigure}
\begin{subfigure}{0.48\textwidth}
\includegraphics[width=1.0\linewidth]{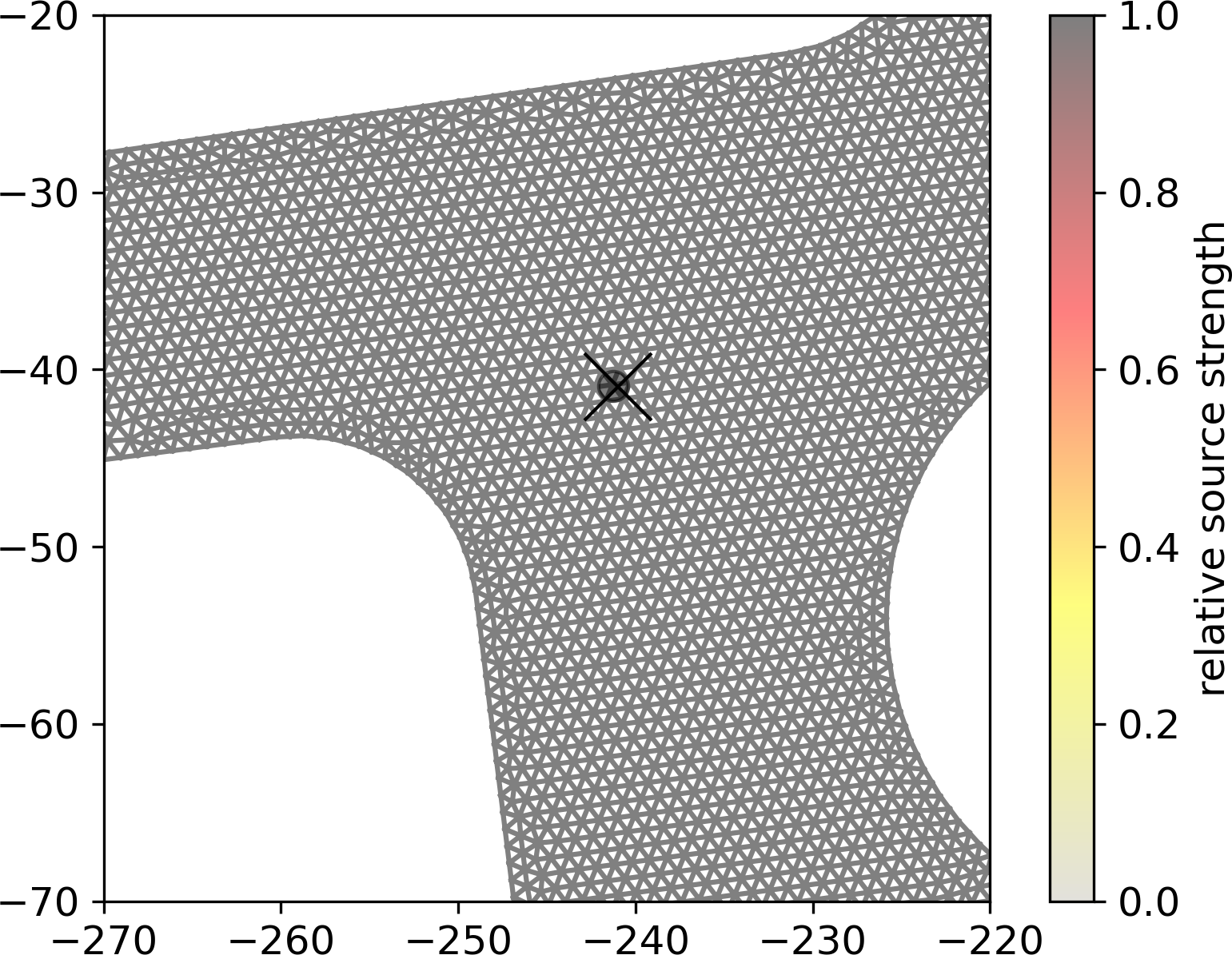}
\end{subfigure}
\caption{Predicted raw sources $\lambda_i\,\delta_{\x^{}_i}$ (left) and derived source $\mu_{\text{post}}$ (red dot) near ground truth source center (black cross, right)}
\label{fig:point_sources}
\end{figure}

\begin{figure}
\centering
\begin{subfigure}{0.48\textwidth}
\includegraphics[width=.9\linewidth]{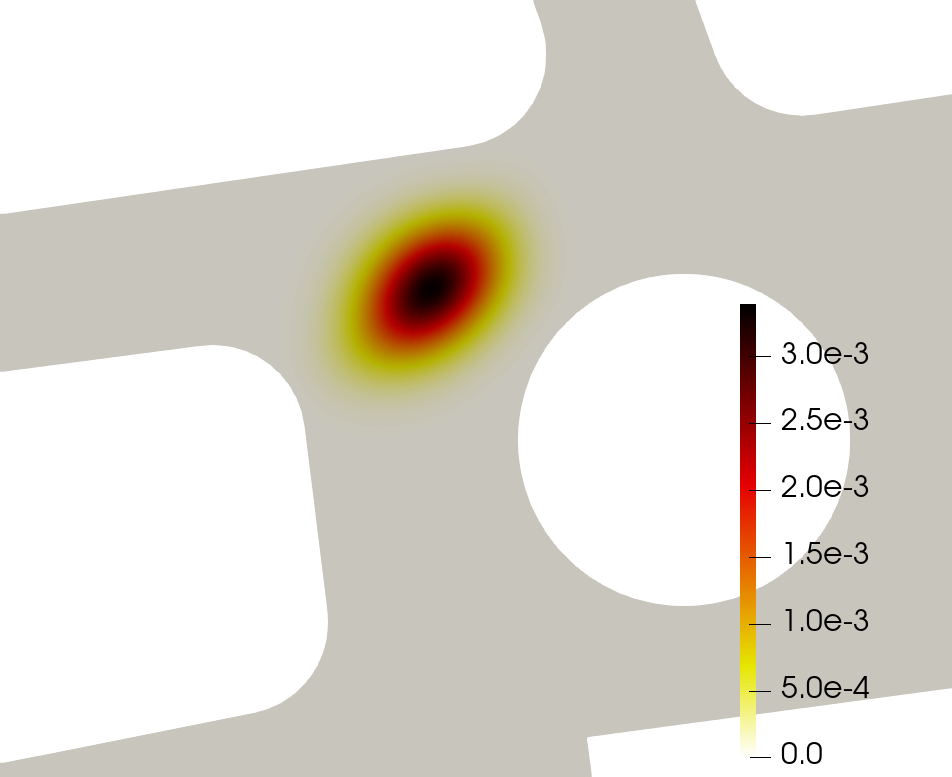}
\end{subfigure}
\begin{subfigure}{0.48\textwidth}
\includegraphics[width=.9\linewidth]{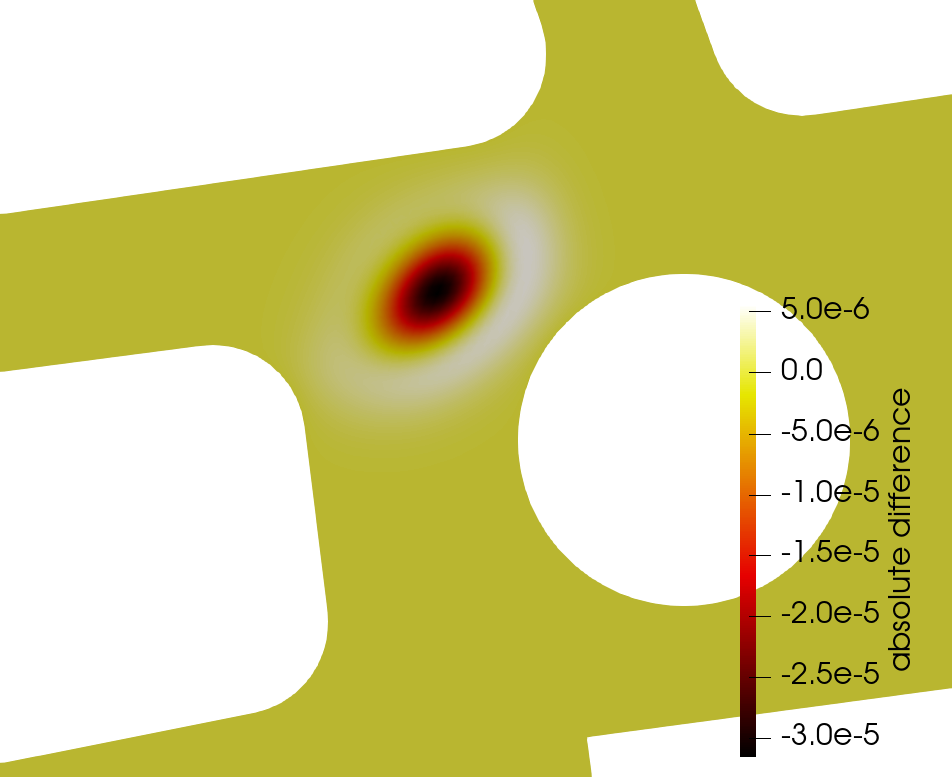}
\end{subfigure}
\caption{Reconstruction of the contaminant concentration at measurement time $\SI{5.0}{\s}$; compare to \autoref{fig:contur} and the absolute difference $u_\text{true}(\SI{5.0}{\s},\cdot)-\hat F(\mu)(\SI{5.0}{\s},\cdot)$ }
\label{fig:tom_reconstruction}
\end{figure}

\begin{figure}
\centering
\begin{subfigure}{0.48\textwidth}
\includegraphics[width=.9\linewidth]{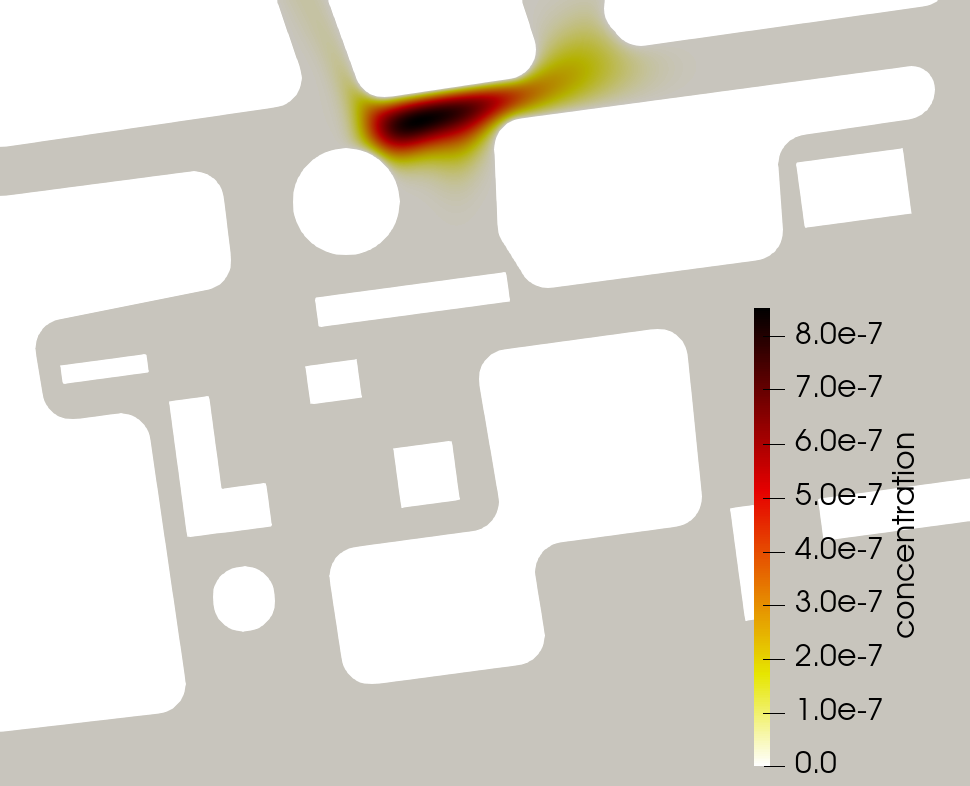}
\end{subfigure}
\begin{subfigure}{0.48\textwidth}
\includegraphics[width=.9\linewidth]{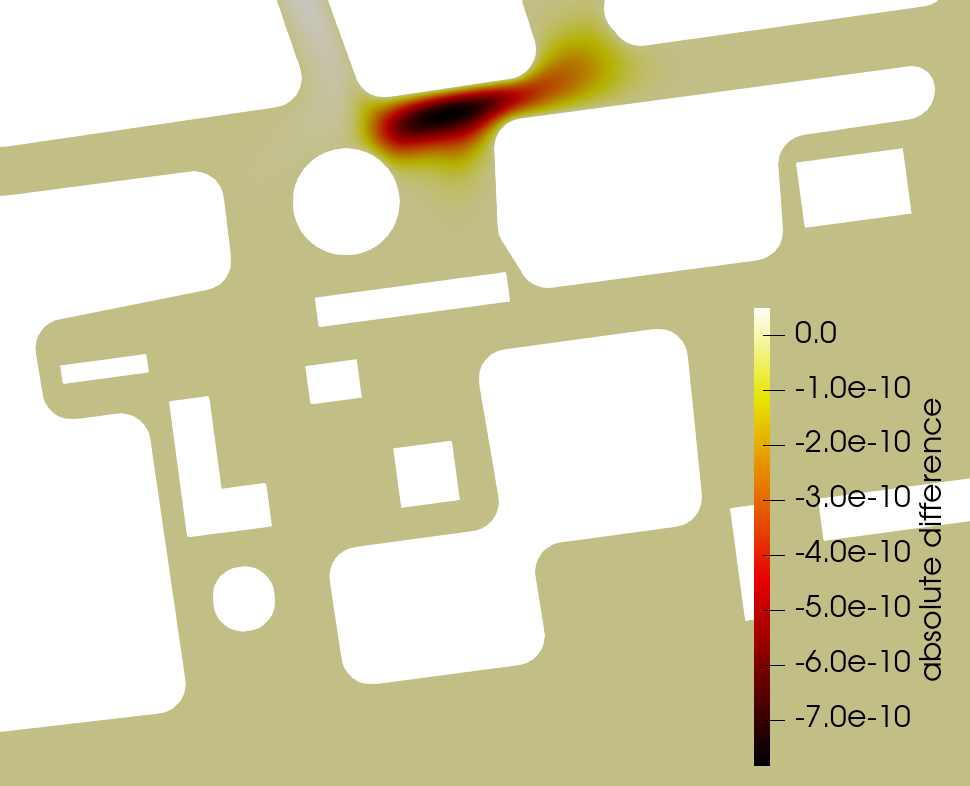}
\end{subfigure}
\caption{Reconstruction of the contaminant concentration at measurement time $\SI{100.0}{\s}$; compare to \autoref{fig:fwd} and the absolute difference $u_\text{true}(\SI{100.0}{\s},\cdot)-\hat F(\mu)(\SI{5.0}{\s},\cdot)$ }
\label{fig:prediction}
\end{figure}
To demonstrate the capabilities of the proposed method, we employ an automated computational pipeline that imports building footprints directly from OpenStreetMap (OSM) and incorporates them as obstacles into the computational domain; see \autoref{fig:domain}. Based on this geometry, locally refined triangular meshes are generated to ensure stable and accurate numerical solutions of both the forward problem (\ref{eq:forward_equation}) and the adjoint problem (\ref{eq:adjoint_equation}). The implementation is based on the software framework \fenics\cite{Baratta.2023}; see also\cite{Bonari.2024} for details of the mesh generation procedure.
\noindent
As a proof of concept, we consider a synthetic test case with artificially generated wind and measurement data. The resulting wind field for a kinematic viscosity of $\nu = \SI{0.25}{\square\metre\per\second}$ is shown in \autoref{fig:domain}. A constant inflow velocity of $\SI{1}{\metre\per\second}$ is prescribed at the southern boundary, while free-slip boundary conditions are imposed elsewhere, cf. \cite{Gjerde.2022}. The corresponding forward simulation is illustrated in \autoref{fig:fwd}. Synthetic measurements are extracted at time $\tobs=\SI{5}{\second}$. To mimic imaging-based sensing techniques, we assume that the measurement device becomes sensitive only for concentrations exceeding the threshold value $c = 10^{-3}$. The corresponding contour is depicted in \autoref{fig:contur}.
\noindent
For the transport problem (\ref{eq:forward_equation}), we choose a diffusion coefficient of $\kappa = \SI{2.0}{\square\metre\per\second}$. The source term is parameterized by the shape function
\begin{equation}
\label{eq:ansatz}
\ansatzSources(\x_s,z)=\min\left\{0.5,\exp\left(\ln(\epsilon)\fracSolidus{\norm{z-\x_s}_2^2}{r^2}\right)\right\},
\end{equation}
with radius $r=\SI{1}{\metre}$, center $\x_s=[-241,-41]$, and threshold parameter $\epsilon = 10^{-3}$. With this setup, the proposed method is able to reconstruct the source location efficiently and with high accuracy. The resulting six Dirac delta distributions $(\lambda_1\delta_{x_1},\dots,\lambda_6\delta_{x_6})$, computed by the PDAP algorithm together with the post-processed source $\mu_{\text{post}}=\lambda_{\text{post}}\delta_{x_{\text{post}}}$, where $\lambda_{\text{post}}=\sum_{i=1}^6 \lambda_i$ and $x_{\text{post}}=\fracSolidus{\sum_{i=1}^6 \lambda_i x_i}{\lambda_{\text{post}}}$, are reported in \autoref{tab:results}. These results were obtained after only sixteen iterations of the PDAP algorithm, requiring a total of sixteen solutions of the forward and adjoint problems \autoref{eq:forward_equation} and \autoref{eq:adjoint_equation}, respectively.
Since both the source location and the source intensity are reconstructed with high accuracy, the complete contaminant distribution at the measurement time can be recovered reliably, enabling subsequent predictions of the future dispersion dynamics. 
For the synthetic test case considered here, the reconstruction error exhibits a spatial distribution similar to that of the concentration field and remains below $1\%$ throughout the domain. This behavior is illustrated in \autoref{fig:tom_reconstruction} and \autoref{fig:prediction}.
\begin{table}
\centering
\begin{tabular}[b]{c c c c c c}
\toprule
 \multicolumn{1}{c}{Source location} & \multicolumn{2}{c}{Reconstruction $\x$}& $\norm{\x-\x_s}_{\Omega}$ & Estimated intensity $\lambda$\\
 \multicolumn{1}{c}{}  & \multicolumn{1}{c}{x} & \multicolumn{1}{c}{y} & \multicolumn{1}{c}{[\si{\m}]} & (truth 1.0)\\
\midrule
$x_{\text{post}},\lambda_{\text{post}}$ &  $-241.2$&$-40.9$ & 0.26 & 1.02 \\
\midrule
$x_{1},\lambda_{1}$ &  $-241.30$&$-40.86$ & 0.33 & 0.9981 \\
$x_{2},\lambda_{2}$ &  $-237.87$&$-43.90$ & 4.27 & 0.0134 \\
$x_{3},\lambda_{3}$ &  $-237.95$&$-39.54$ & 3.39 & 0.0110 \\
$x_{4},\lambda_{4}$ &  $-240.61$&$-46.01$ & 5.03 & 0.0042 \\
$x_{5},\lambda_{5}$ &  $-241.99$&$-35.71$ & 5.38 & 0.0006 \\
$x_{6},\lambda_{6}$ &  $-247.40$&$-44.30$ & 7.20 & 0.0007 \\

\bottomrule
\end{tabular}
\caption{Reconstruction of the source location  $\x_s=[-241,-41]$ \tcr{at} the test site, cf. \autoref{fig:point_sources}}
\label{tab:results}
\end{table}
\section{Summary and Further Application}\label{sec:outlook}
\begin{figure}
\begin{subfigure}{0.70\textwidth}
\centering
\includegraphics[width=1.0\linewidth]{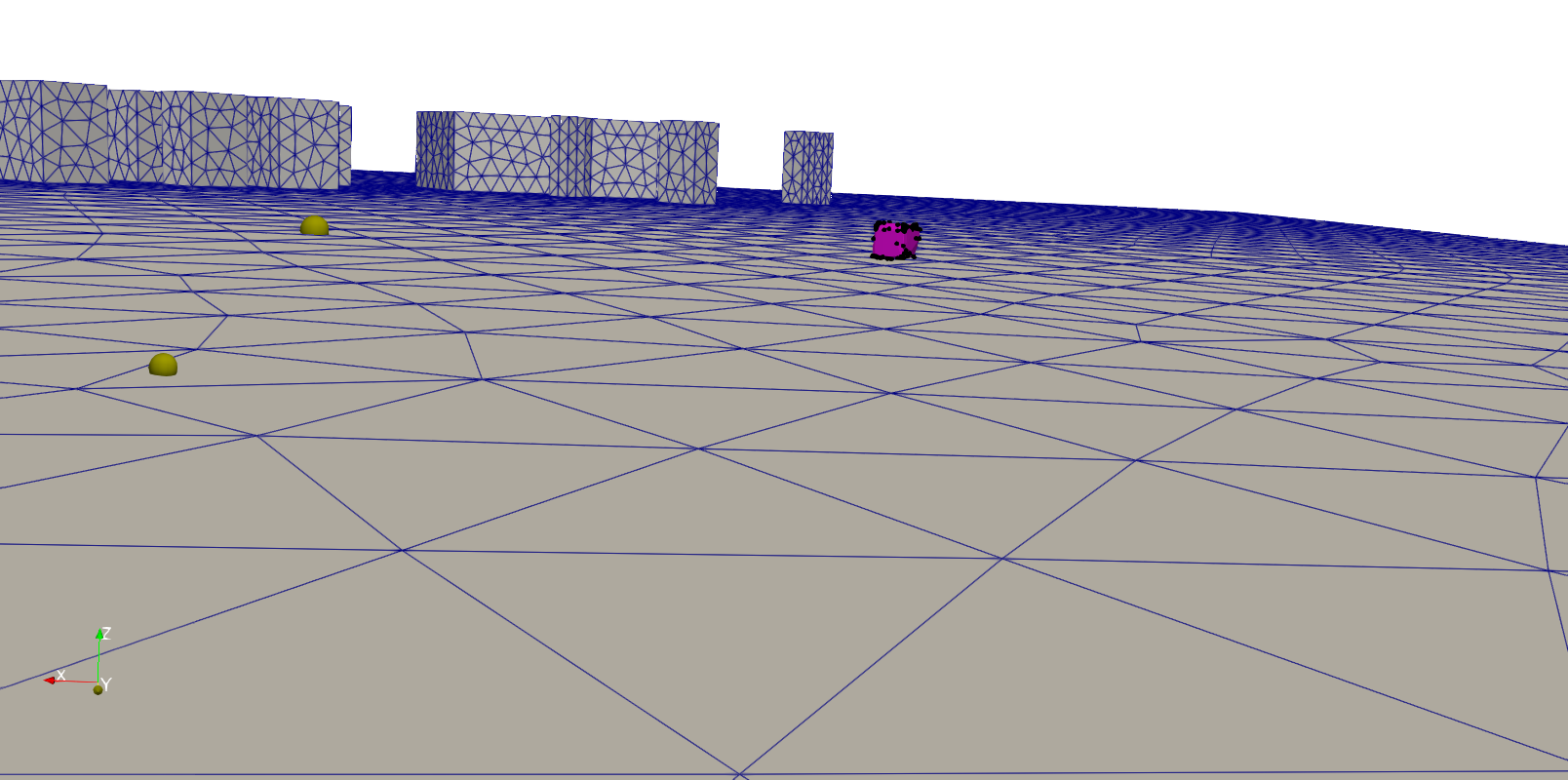}
\end{subfigure}
\begin{subfigure}{0.29\textwidth}
\centering
\includegraphics[width=1.0\linewidth]{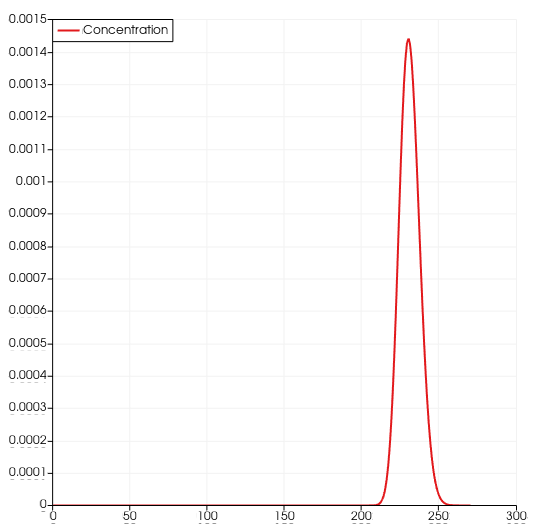}
\end{subfigure}
\caption{Shows the recorded values (black dots) of two \tcr{FPA-FTIR hyperspectral} imaging devices (green-yellow dots) and the extracted contour (pink) in a three-dimensional environment (left) and a synthetic signal representing the information gained from a lidar measurement in the modeled setup (right)}
\label{fig:real_measurement}
\end{figure}
This paper extends a novel and computationally efficient framework for the real-time \tcr{localization} of contaminant sources in complex environments using tomographic \tcr{standoff detection}. In particular, the proposed approach can incorporate high-resolution \tcr{hyperspectral} imaging data acquired by \tcr{an FPA-}FTIR system (so far simulated) to localize and quantify airborne contaminant releases.
A key advantage of the presented imaging technique is its mobile applicability combined with its ability to provide high spatial data density. A real (physical) measurement of a contaminant cloud captured with \tcr{two FPA-FTIR hyperspectral} imaging systems is shown in \autoref{fig:real_measurement}. Compared to isolated point measurements, \tcr{FPA-FTIR hyperspectral FPA-FTIR hyperspectral imaging systems provide} significantly more information about the spatial distribution of the contaminant field. These properties are particularly important for early warning and rapid emergency response in hazardous substance release scenarios.

The computational efficiency of the present formulation enables the identification of sources on large spatial scales. Even for high-resolution measurement data, the proposed approach requires only a comparatively small number of PDE solves. This efficiency is essential for extending the approach to real measurement data in complex three-dimensional environments. To further accelerate the computations, a reduced-order modeling approach for the advection-diffusion problem was investigated in\cite{MattuschkaGoal}. Using singular value decomposition (SVD), a reduced-order model can be constructed that enables highly efficient evaluation of both the forward and adjoint problems and may additionally be exploited for tasks such as optimal sensor placement. In a comparable setting, the reduced model required only $6.25\,\mathrm{ms}$ per forward or adjoint solve. However, such a reduction approach cannot be transferred directly to the present application in a meaningful way, since the boundary curve $\Gamma_{\tobs}[u_{\text{true}}]$ depends on the current measurement and therefore cannot be precomputed offline. Consequently, alternative strategies, e.g., machine learning, are required to achieve real-time capabilities for the proposed method.

Technical measurement systems of this type are inherently affected by systematic noise. In the present work, sensor noise has not yet been incorporated into the inverse problem. A possible formulation based on noisy level-set observations is given by
\begin{equation*}
\Gamma^{\text{obs}}=\left\{x \mid u(x) + \eta(x) = c\right\},\qquad\eta(x) \sim \mathcal{N}(0,\sigma^2).
\end{equation*}
To calibrate and validate such a noise model, data from existing measurement campaigns may be analyzed.
In addition, future work will focus on combining multiple measurement systems in order to obtain reliable data under realistic environmental conditions. In practical applications, not only the initial condition but also additional unknown model parameters must be identified through data assimilation. Since detailed flow simulations are often infeasible due to limited data availability and the complexity of the environment, additional parameters related to turbulence modeling must be estimated, which may require a simultaneous calibration of the overall system. This leads to the fact that data fusion of different types of measurements, e.g., combining \tcr{standoff FPA-FTIR} and \tcr{portable polarization lidar (light imaging, detection and ranging) system}~\cite{Zajonz.2026}, allows for a significant gain in information and thus lead to more precise physically grounded statements for real-time situational awareness in hazardous substance release scenarios. \tcr{An extracted synthetic example of the received power signal from the portable polarization lidar system (\autoref{fig:hi90}, left, white line) is shown in \autoref{fig:real_measurement}, where attenuation effects, i.e., the attenuation of the emitted and backscattered laser signal within the plume due to absorption and scattering, are neglected. Future work will focus on incorporating these effects in a physically consistent manner, including an adequate noise model.}
\section{Acknowledgements}
 AP gratefully acknowledges the funding by dtec.bw - Digitalization and Technology Research Center of the Bundeswehr (project RISK.twin). dtec.bw is funded by the European Union - NextGenerationEU. 

\printbibliography

\end{document}